\documentclass[12pt]{amsart}

\usepackage{amsfonts,amssymb,amsthm}

\newtheorem{thm}{Theorem}
\newtheorem{cor}[thm]{Corollary}

\newtheorem{prop}[thm]{Proposition}

\newtheorem{ex}[thm]{Example}

\newcommand{\R}{\mathbb{R}}
\newcommand{\E}{\mathbb{E}}
\renewcommand{\P}{\mathbb{P}}
\newcommand{\N}{\mathbb{N}}
\newcommand{\K}{\mathcal{K}}
\DeclareMathOperator{\vol}{vol}
\DeclareMathOperator{\conv}{conv}
\newcommand{\bp}{\begin{pmatrix}}
\newcommand{\ep}{\end{pmatrix}}
\newcommand{\Par}{\mathcal{P}}
\newcommand{\Ell}{\mathcal{E}}
\renewcommand{\epsilon}{\varepsilon}

\numberwithin{thm}{section}

\title{On symmetric versions of Sylvester's problem}

\author[M. Meckes]{Mark W. Meckes}

\address{Department of Mathematics, Case Western Reserve University,
    Cleveland, Ohio 44106.}

\email{mwm2@po.cwru.edu}


\thanks{Partially supported by a grant from the National Science Foundation.}

\begin{document}

\begin{abstract}
We consider moments of the normalized volume of a symmetric or
nonsymmetric random polytope in a fixed symmetric convex body. 
We investigate for which bodies these moments are extremized, and
calculate exact values in some of the extreme cases. We show
that these moments are maximized among planar convex bodies by
parallelograms.
\end{abstract}

\maketitle

\section{Introduction} \label{S:Intro}

Sylvester's four point problem asks for the probability that the convex
hull of four random points, chosen independently and uniformly from a
convex body $K \subset \R^2$, is a quadrilateral, and in particular, for
which convex bodies $K$ this probability is extremal. This is equivalent
to asking what the expected area of the convex hull of three random points
in $K$ is, and for which bodies this expectation is extremal. This problem
was solved by Blaschke \cite{Blaschke1, Blaschke2}, who showed that the
expected area achieves its maximum exactly when $K$ is a triangle and
achieves its minimum exactly when $K$ is an ellipse.

Since then, various authors have considered several extensions of this
problem. Many of these are special cases of the following general problem.
We write $\K^n$ for the set of all convex bodies in $\R^n$, that is, all
compact convex sets with interior points. Let $K \in \K^n$ and $N \ge
n+1$. Let $x_1, x_2, \ldots, x_N$ be independent random points distributed
uniformly in $K$. We define
$$U_{K, N} = \frac{\vol_n(\conv \{ x_1, x_2, \ldots, x_N \})}
    {\vol_n(K)};$$
thus the random variable $U_{K,N}$ is the normalized volume of a random
polytope in $K$. Note that the distribution of $U_{K, N}$ is an affine
invariant of $K$. The generalized Sylvester's problem asks, for each $n
\ge 2$, $N \ge n+1$, and $p \ge 1$, for which $K \in \K^n$ does the $p$th
moment $\E U_{K, N}^p$ achieve its extremal values? It should be noted at
this point that a compactness argument guarantees that such extremal
bodies do exist; see \cite{Groemer1, DL, Giann, CCG}.

Groemer \cite{Groemer1, Groemer2} showed that, for each such $n$, $N$, and
$p$, $\E U_{K, N}^p$ is minimized exactly when $K$ is an ellipsoid. Dalla
and Larman \cite{DL} showed for $n = 2$ that $\E U_{K, N}$ is maximized,
for each $N \ge 3$, when $K$ is a triangle; and Campi, Colesanti, and
Gronchi \cite{CCG} extended this to $\E U_{K, N}^p$ for all $p \ge 1$.
Giannopoulos \cite{Giann} showed that for $n = 2$, $\E U_{K, N}$ is
maximized only if $K$ is a triangle. Very little is known about maximizing
bodies when $n \ge 3$. It is widely conjectured that $\E U_{K, N}^p$
should achieve its maximum exactly when $K$ is a simplex, but there are
only partial results in this direction \cite{DL, CCG}. As noted explicitly
in \cite{Meckes} (but see also \cite[Proposition 5.6]{MP}), this would in
particular imply the well-known hyperplane conjecture \cite{MP, GM}.

In this paper, we consider two ``symmetric" variants of this generalized
Sylvester's problem. We write $\K_s^n$ for the set of symmetric convex
bodies in $\R^n$; that is, all $K \in \K^n$ such that $K = -K$. The first
variant asks, for which $K \in \K_s^n$ does $\E U_{K, N}^p$ achieve its
extremal values? The second variant asks the same question when the random
polytope, as well as the fixed body $K$, is symmetric. More precisely, for
$K \in \K_s^n$ and $N \ge n$, we again let $x_1, x_2, \ldots, x_N$ be
independent random points distributed uniformly in $K$. We define
$$V_{K, N} = \frac{\vol_n(\conv \{ \pm x_1, \pm x_2, \ldots, \pm x_N \})}
    {\vol_n(K)}.$$
The distribution of $V_{K, N}$ is a linear invariant of $K$. We now ask,
for each $N \ge n$ and $p \ge 1$, for which $K \in \K_s^n$ does $\E V_{K,
N}^p$ achieve its extremal values?

The first goal of this paper is to bring the level of knowledge about
these symmetric versions of Sylvester's problem to a level close to that
for the nonsymmetric case. Since $\E U_{K,N}^p$ is already known to be
minimized over all $K \in \K^n$ exactly when $K$ is an ellipsoid, it is in
particular minimized over all $K \in \K_s^n$ exactly when $K$ is an
ellipsoid. Furthermore, it was noted in \cite{Meckes} that Groemer's proof
also shows that $\E V_{K,N}^p$ is minimized over all $K \in \K_s^n$
exactly when $K$ is an ellipsoid. Thus in this paper we will deal with the
question of which $K \in \K_s^n$ maximize $\E U_{K, N}^p$ and $\E V_{K,
N}^p$. We show in Theorem \ref{thm:max} that when $n  = 2$, each maximum
is achieved when $K$ is a parallelogram. Our main tools, which we
introduce in Section \ref{S:SRS}, are symmetric adaptations of tools
developed by Campi, Colesanti, and Gronchi \cite{CCG} to study the
nonsymmetric generalized Sylvester's problem. Following \cite{CCG}, we
derive some partial results for general $n$, which in particular support
the conjecture that the maximizing symmetric convex bodies should be
either parallelotopes or crosspolytopes, or bodies built from these.

The second goal of this paper is to obtain information about the extremal
values of $\E V_{K,N}^p$, which we do in Section \ref{S:Calculations}.
When $n = 2$, we derive the exact distributions of the random variables
$V_{\Par, 2}$ and $V_{\Ell^2, 2}$, where $\Par$ denotes a parallelogram
and $\Ell^2$ denotes an ellipse; and we calculate $\E V_{\Par, N}$ and $\E
V_{\Ell^2, N}$ for all $N \ge 2$. We also calculate $\E V_{\Ell^3, N}$ for
all $N \ge 3$, where $\Ell^3$ denotes an ellipsoid in $\R^3$. The
corresponding extremal values of $\E U_{K,N}^p$ are already available in
the literature.

\section{RS- and SRS-decomposability} \label{S:SRS}

In this section, we recall the notions of RS-movements and
RS-decomposability of a convex body, which were introduced by Campi,
Colesanti, and Gronchi in \cite{CCG}, and introduce complementary notions
for symmetric convex bodies. These tools will be used to address the
problem of identification of maximizers of $\E V_{K, N}^p$ and $\E U_{K,
N}^p$ for $K \in \K_s^n$.

We first recall the notion of a linear parameter system, due to Rogers and
Shephard \cite{RS}. For $n \ge 2$, let $K \in \K^n$, $\alpha : K \to \R$,
and $v \in \R^n \setminus \{ 0 \}$. Then for each $t$ in some interval in
$\R$, we set
$$ K_t = \conv \{ x + t \alpha(x) v : x \in K \} .$$
The family of sets $K_t$ is called a {\em linear parameter system} with
speed function $\alpha$. The most important property of linear parameter
systems is the following, proved in \cite{RS}.

\begin{thm}[Rogers - Shephard] \label{thm:RS}
Let $K_t$ be a linear parameter system for $K \in \K^n$.  Then $\vol_n
(K_t)$ is convex as a function of $t$.
\end{thm}

As in \cite{CCG}, the interest here is in the case in which the speed
function is constant on each chord of $K$ which is parallel to $v$. Let
$\pi_v : \R^n \to v^\perp$ denote orthogonal projection, and let $\beta :
\pi_v (K) \to \R$. In the terminology of \cite{CCG}, a family of sets
\begin{equation} \label{eq:RS-move1}
K_t = \{ x + t \beta(\pi_v(x)) v : x \in K \},
\end{equation}
for all $t$ in some interval containing 0, is called an {\em RS-movement}
of $K$ if $K_t$ is convex for each allowed $t$. A more convenient way to
describe $K_t$ is the following. Given $v \in \R^n$, there exist functions
$f_v, g_v : \pi_v (K) \to \R$ with $f_v$ convex and $g_v$ concave such
that
\begin{equation} \label{eq:convex}
K = \{ x + r v : x \in \pi_v (K), \, f_v(x) \le r \le g_v(x) \}.
\end{equation}
Then if $\beta : \pi_v(K) \to \R$ is given, \eqref{eq:RS-move1} is
equivalent to
\begin{equation} \label{eq:RS-move2}
K_t = \{ x + r v : x \in \pi_v (K), \,
    f_v(x) + t \beta(x) \le r \le g_v(x) + t \beta(x) \}.
\end{equation}
Note that a necessary and sufficient condition for $\beta : \pi_v (K) \to
\R$ to define an RS-movement is that $f_v + t \beta$ is convex and $g_v +
t \beta$ is concave for each allowed $t$.

It is noted in \cite{CCG} that if $\beta$ is any affine function on $\pi_v
(K)$, then $K_t$ as defined by \eqref{eq:RS-move1} is an RS-movement of
$K$ such that each $K_t$ is an affine image of $K$. Moreover, Steiner
symmetrization is related to a particular RS-movement as follows. If
$\beta = - (f_v + g_v)$, we obtain an RS-movement of $K$ such that $K_1$
is the reflection of $K$ with respect to $v^\perp$, and $K_{1/2}$ is the
Steiner symmetrization of $K$ with respect to $v^\perp$.

Now let $K \in \K_s^n$. We say that an RS-movement $K_t$ of $K$ is an {\em
SRS-movement} if the speed function $\beta : \pi_v (K) \to \R$ is odd,
that is, if $\beta(-x) = -\beta(x)$. Note that this is precisely the
condition which ensures that $K_t \in \K_s^n$ for each $t \in [a, b]$.
Note that if $\beta$ is any linear function on $v^\perp$, then $K_t$ as
defined by \eqref{eq:RS-move1} is an SRS-movement of $K$ such that each
$K_t$ is a linear image of $K$. Furthermore, if $K$ is symmetric, then for
any $v \in \R^n$, the functions $f_v, g_v$ in \eqref{eq:convex} satisfy
$g_v(-x) = -f_v(x)$. For example, the RS-movement with speed function
$\beta = -(f_v + g_v)$, which gives rise to reflection and Steiner
symmetrization with respect to $v^\perp$, is an SRS-movement of $K$.

Following \cite{CCG}, we say that $K \in \K^n$ is {\em RS-decomposable} if
there exists an RS-movement $K_t$, $t \in (-\epsilon, \epsilon)$ for some
$\epsilon > 0$, such that $K_0 = K$ and such that the speed function is
not affine. $K \in \K^n$ is called {\em RS-indecomposable} if it is not
RS-decomposable. In analogy, we say that $K \in \K_s^n$ is {\em
SRS-decomposable} if there exists an SRS-movement $K_t$, $t \in
(-\epsilon, \epsilon)$ for some $\epsilon > 0$, such that $K_0 = K$ and
such that the speed function is not linear. $K \in \K_s^n$ is called {\em
SRS-indecomposable} if it in not SRS-decomposable.

We remark at this point that to avoid ambiguity, we maintain a strict
distinction between affine and linear functions, even in one dimension, so
that a linear function $f:\R \to \R$ is required to satisfy $f(0) = 0$.

\begin{ex} \label{ex:square}
A symmetric parallelogram $\Par \in \K_s^2$ is SRS-indecomposable.
\end{ex}
\begin{proof}
We may identify $v^\perp$ with $\R$. Then there exist $0 \le a \le b$ such
that $\pi_v (\Par) = [-b, b]$ and such that the functions $f_v, g_v$ as in
\eqref{eq:convex} are affine on each of the intervals $[-b, -a]$, $[-a,
a]$, and $[a, b]$. Moreover, one of $f_v, g_v$ is affine on $[-a, b]$ and
the other is affine on $[-b, a]$. Assume without loss of generality that
$f_v$ is affine on $[-a, b]$. In order for $f_v + t \beta$ to be convex
for both positive and negative values of $t$, $\beta$ must be linear on
$[-a, b]$. Since $\beta$ is odd, this implies that $\beta$ is linear on
$\pi_v (\Par)$.
\end{proof}

If $K = \conv (K' \cup \{ x, -x \} )$, where $K'$ is a symmetric convex
body in a hyperplane $H$ and $x \notin H$, then we call $K$ a {\em double
cone} with base $K'$.

\begin{prop} \label{prop:cyl-cone}
Let $K \in \K_s^n$ be either a symmetric cylinder or a double cone. Then
$K$ is SRS-decomposable if and only if its base is.
\end{prop}

The proof of this is an almost verbatim repetition of \cite[Example
2.6]{CCG}, which shows that a cylinder or cone is RS-decomposable if and
only if its base is. Combined with Example \ref{ex:square}, Proposition
\ref{prop:cyl-cone} implies the following.

\begin{cor} \label{cor:cube-crosspoly}
Any symmetric parallelotope or crosspolytope in $\R^n$, $n \ge 2$, is
SRS-indecom\-posable.
\end{cor}

We note that the results of \cite{CCG} imply that every simplex is
RS-indecomposable, whereas every parallelotope is RS-decomposable.

The proof of \cite[Theorem 3.3]{CCG} also yields the following.

\begin{prop} \label{prop:curvature}
Let $K \in \K_s^n$ be such that $\partial K$ has a nonempty open subset of
class $C^2$ on which all the principal curvatures are positive. Then $K$
is SRS-decomposable.
\end{prop}

The main technical result of \cite{CCG} is the following.

\begin{prop}[Campi-Colesanti-Gronchi]
Let $K_t$ be an RS-movement of $K \in \K^n$. Then $\E U_{K_t, N}^p$ is a
convex function of $t$, for every $p \ge 1$ and $N \ge n + 1$.
Furthermore, $\E U_{K_t, N}^p$ is strictly convex if and only if the speed
function is not affine.
\end{prop}

Theorem \ref{thm:RS} is the main tool used to prove this. With minor
modifications, the same proof yields the following.

\begin{prop} \label{prop:SRS}
Let $K_t$ be an SRS-movement of $K \in \K_s^n$. Then $\E V_{K_t, N}^p$ is
a convex function of $t$, for every $p \ge 1$ and $N \ge n$. Furthermore,
$\E V_{K_t, N}^p$ is strictly convex if and only if the speed function is
not linear.
\end{prop}

As an immediate consequence of Proposition \ref{prop:SRS} and the
definition of SRS-decomposability, we have the following.

\begin{cor} \label{cor:SRS}
Let $p \ge 1$ and $N \ge n$. If $K$ maximizes $\E U_{K, N + 1}^p$ or $\E
V_{K, N}^p$ for all $K \in \K_s^n$, then $K$ is SRS-indecomposable.
\end{cor}

Corollary \ref{cor:SRS} and Proposition \ref{prop:curvature} suggest (but
do not imply) that the maximizers of $\E U_{K, N + 1}^p$ and $\E V_{K,
N}^p$ are polytopes. Furthermore, Corollary \ref{cor:cube-crosspoly} shows
that the present method will not rule out the obvious candidates.

\begin{thm} \label{thm:max}
For any $K \in \K_s^2$, $p \ge 1$, and $N \ge 2$,
\begin{align*}
\E U_{K, N + 1}^p & \le \E U_{\Par, N + 1}^p, \\
\E V_{K, N}^p     & \le \E V_{\Par, N}^p,
\end{align*}
with strict inequality in both of the above if $K$ is a symmetric polygon
with more than 4 vertices.
\end{thm}

\begin{proof}
Suppose that $K$ is a symmetric polygon with vertices $\pm P_1, \pm P_2,
\ldots, \pm P_m$, $m \ge 3$, ordered so that $P_i, P_{i+1}$ are adjacent
for each $i = 1, 2, \ldots, m - 1$. Then $P_2, -P_m$ are the vertices
adjacent to $P_1$. We set
$$ K_t = \conv \bigl\{ \pm \bigl(P_1 + t (P_2 + P_m)\bigr), \pm P_2, \pm P_3,
    \ldots, \pm P_m \bigr\}.$$
There exist an $\epsilon_1 > 0$ such that $P_1 + \epsilon_1 (P_2 + P_m)$
lies on the line through $P_2$ and $P_3$, and an $\epsilon_2 > 0$ such
that $P_1 - \epsilon_2 (P_2 + P_m)$ lies on the line through $-P_m$ and
$-P_{m-1}$. Then $K_t$, $t \in [-\epsilon_2, \epsilon_1]$ is an
SRS-movement such that $K_{\epsilon_1}$ and $K_{-\epsilon_2}$ have $2 (m -
1)$ vertices. Furthermore, this SRS-movement fixes
$$ \conv \{ \pm P_2, \pm P_3, \ldots, \pm P_m \}, $$
and therefore the corresponding speed function is not linear. Thus
Proposition \ref{prop:SRS} implies
\begin{align*}
\E U_{K, N+1}^p & < \max \{ \E U_{K_{\epsilon_1}, N + 1}^p,
    \E U_{K_{-\epsilon_2}, N + 1}^p \}, \\
\E V_{K, N}^p & < \max \{ \E V_{K_{\epsilon_1}, N}^p,
    \E V_{K_{-\epsilon_2}, N}^p \}.
\end{align*}
Iterating this argument, we obtain $\E U_{K, N + 1}^p < \E U_{\Par, N +
1}^p$ and $\E V_{K, N}^p < \E V_{\Par, N}^p$.

The statement for a general $K \in \K_s^2$ now follows by the continuity
of $\E U_{K, N + 1}^p$ and $\E V_{K, N}^p$ as functions of $K$.
\end{proof}

\begin{cor} \label{cor:isotropic}
For $K \in \K^n$, let $L_K$ denote the isotropic constant of $K$.
\begin{enumerate}
\item For any $K \in \K_s^2$, $L_K \le L_{\Par} = (12)^{-1/2}$.
\item For any $K \in \K^2$ with centroid at the origin, $L_K \le L_\Delta
= (108)^{-1/4}$, where $\Delta$ denotes a triangle with centroid at the
origin.
\end{enumerate}
\end{cor}

\begin{proof}
The first claim follows directly from Theorem \ref{thm:max} and the
formula
$$ L_K^2 = \frac{1}{4} \left(n! \E V_{K, n}^2\right)^{1/n} $$
for any $K \in \K_s^n$ (see \cite{Ball, Meckes}). The second claim follows
from the fact that $\E U_{K,N}^p \le \E U_{\Delta,N}^p$ for $N \ge 3$ and
$p \ge 1$ \cite{CCG} and the formula
$$ L_K^2 = \left(\frac{n! \E U_{K, n + 1}^2}{n + 1} \right)^{1/n}$$
for any $K \in \K^n$ with centroid at the origin, due (essentially) to
Kingman \cite{Kingman}.
\end{proof}

Schmuckenschl\"ager proved \cite{Schmuck} that for all $n$, $L_{B_p^n} \le
L_{B_\infty^n} = (12)^{-1/2}$ for all $1 \le p \le \infty$, where $B_p^n$
is the unit ball of $\ell_p^n = (\R^n, \| \cdot \|_p)$. This fact and
Corollary \ref{cor:isotropic} support the conjecture $L_K \le (12)^{-1/2}$
for all $K \in \K_s^n$, $n \in \N$. This may be considered an isometric
form of the hyperplane conjecture.

\section{Calculations for parallelograms and ellipsoids}
\label{S:Calculations}

In this section we calculate some extremal values of $\E V_{K,N}^p$. In
Section \ref{S:Densities}, we derive the exact distributions of $V_{K,2}$
when $K$ is either a parallelogram or an ellipse, making essential use of
the symmetries of those bodies. In Sections \ref{S:Ell-Par} and
\ref{S:Ellipsoid}, we derive general formulas for $\E V_{K,N}$ for $N \ge
n$ and $n = 2, 3$ respectively. When $n = 2$ we use these to derive simple
expressions in the cases of parallelograms and ellipses; when $n = 3$ we
derive an expression for ellipsoids. We also indicate where the
corresponding values of $\E U_{K, N}^p$ may be found in the literature.

We remark that if $\Ell^n$ denotes an ellipsoid in $\R^n$, $\E V_{\Ell^n,
n}^p$ was computed for $n \in \N$ and $p > 0$ by the author in
\cite{Meckes}, and $\E U_{\Ell^n, n+1}^p$ was computed for $n, p \in \N$
by Miles in \cite{Miles}.

\subsection{Densities when $n = N = 2$} \label{S:Densities}

\begin{prop}
$V_{\Par,2}$ has density
$$ I_{[0, 1]}(t) \int_{2t - 1}^1 (\log |s|) (\log |2t - s|) ds.$$
\end{prop}
\begin{proof}
We may assume that $\Par$ is the square $[-1,1]^2$. Since the symmetric
convex hull of two points $x, y \in \R^2$ has area $2 |x_1 y_2 - x_2
y_1|$, $V_{\Par,2}$ has the same distribution as $\frac{1}{2} |X_1 X_2 -
X_3 X_4|$, where $X_i$, $1 \le i \le 4$, are independent random variables
uniformly distributed in $[-1, 1]$. By symmetry, $V_{\Par,2}$ also has the
same distribution as $\frac{1}{2} |X_1 X_2 + X_3 X_4|$. We begin by
calculating the distribution of $X_1 X_2$. First note that $X_1 X_2$ is
symmetric. Now, for $t > 0$,
$$\P [ X_1 X_2 \le t] = \frac{1}{2} (1 + a_t),$$
where $a_t$ is the area of $\{ (x_1, x_2) \in [0, 1]^2 : x_1 x_2 \le t
\}$. By elementary integration, we obtain $a_t = t (1 - \log t)$ for $0 <
t \le 1$, and $a_t = 1$ for $t > 1$. From this we obtain that $X_1 X_2$
has density
$$ \frac{d}{dt} \P[X_1 X_2 \le t] = -\frac{1}{2} \log |t|,$$
supported on $[-1, 1]$. The distribution of $X_1 X_2 + X_3 X_4$ is then
the convolution of this distribution with itself, so its density is
$$ f(t) = \frac{1}{4} \int_{\max \{ -1, t-1\}}^{\min \{1, t+1\}}
    (\log |s|) (\log |t-s|) ds, $$
supported on $[-2, 2]$. Finally, $\frac{1}{2} |X_1 X_2 - X_3 X_4|$ has
density
$$ 4 f(2t) = \int_{2t-1}^1 (\log |s|) (\log |2t-s|) ds,$$
supported on $[0, 1]$.
\end{proof}

\begin{prop}
$V_{\Ell^2,2}$ has density
$$\pi t I_{[0, \frac{2}{\pi}]}(t) \int_{\frac{\pi}{2} t}^1 s^{-2}
    \sqrt{1 - s^2} ds.$$
\end{prop}
\begin{proof}
We may assume that $\Ell^2$ is the unit disc. By the rotational invariance
of the uniform measure on $\Ell^2$, $V_{\Ell^2,2}$ has the same
distribution as the $\frac{1}{\pi}$ times the area of the symmetric convex
hull of two independent random points, one uniformly distributed in
$\Ell^2$, the other distributed in the interval $[0,1]$ on the $x$-axis
with density $2 t$. Note that since one of the random points lies on the
$x$-axis, the area of their symmetric convex hull depends only on the
absolute value of the $y$-coordinate of the other point, which is
distributed in $[0,1]$ with density $\frac{2}{\pi} \sqrt{1 - t^2}$.
Therefore $V_{\Ell^2,2}$ has the same distribution as $\frac{2}{\pi} XY$,
where $X$ and $Y$ are independent random variables in $[0,1]$ such that
$X$ has density $2 t$ and $Y$ has density $\frac{2}{\pi} \sqrt{1- t^2}$.
$V_{\Ell^2,2}$ then has density supported on $[0,\frac{2}{\pi}]$ given by
\begin{equation*} \begin{split}
\frac{d}{dt} \P \left[XY \le \frac{\pi}{2} t \right] & =
    \frac{d}{dt} \int_0^1 \P \left[X \le \frac{\pi t}{2s} \right]
    \frac{2}{\pi} \sqrt{1 - s^2} ds \\
& = \int_0^1 I_{[0,1]} \left(\frac{\pi t}{2 s}\right) \pi t s^{-2}
    \sqrt{1 - s^2} ds.
\end{split} \end{equation*}
\end{proof}

Exact densities of $U_{K, N}$ have not been derived; however, $\E
U_{\Delta, 3}^p$ and $\E U_{\Par, 3}^p$ were calculated for all $p \in \N$
by Reed \cite{Reed}. The values of $\E U_{\Ell^2, 3}^p$ for $p \in \N$ are
a special case of the above mentioned result of Miles \cite{Miles}.

\subsection{Expected area in an ellipse or parallelogram} \label{S:Ell-Par}

In this and the next section we derive general formulas for $\E V_{K, n}$
when $n = 2, 3$. The derivations make use of standard arguments for
geometric probability, adapted for the symmetric case; see for example the
papers of R\'enyi and Sulanke \cite{RS1, RS2} and Buchta and Reitzner
\cite{BR} for related formulas derived using similar ideas. Our
derivations follow the outline of Buchta and Reitzner's proof of a
nonsymmetric analogue of Proposition \ref{n=3} below.

Let $K \in \K_s^2$. For $r\ge 0$ and $0\le \theta < 2\pi$, let
$$\ell(r, \theta) = \vol_1 (\{(x, y) : (x, y) \cdot
    (\cos \theta, \sin \theta) = r \} \cap K),$$
\begin{equation*} \begin{split}
A(r, \theta) & = \vol_2 (\{(x, y) : |(x, y) \cdot
    (\cos \theta, \sin \theta)| \le r \} \cap K) \\
& = 2\int_0^r \ell(s, \theta) ds.
\end{split} \end{equation*}

\begin{prop}\label{n=2}
Let $K \in \K_s^2$ and $N \ge 2$. If $|K| = 1$, then
$$\E V_{K, N} = 1 - \frac{N}{3} \int_0^{2\pi} \int_0^\infty
    A(r, \theta)^{N - 1} \ell(r, \theta)^3 dr d\theta.$$
\end{prop}
\begin{proof}
We consider a random convex polygon $\Pi_{N + 1}$ which is the symmetric
convex hull of $N + 1$ independent random points distributed uniformly in
$K$. Each of these random points is a vertex of $\Pi_{N + 1}$ iff it is
not contained in the symmetric convex hull of the other $N$ random points,
therefore it is a vertex with probability $1 - \E V_{K, N}$. Each of the
random points is also a vertex iff its antipode also is. Therefore the
expected number $v_{N + 1}$ of vertices of $\Pi_{N + 1}$ is
$$v_{N + 1} = 2 (N + 1)(1 - \E V_{K, N}),$$
and thus
$$\E V_{K, N} = 1 - \frac{v_{N+1}}{2 (N + 1)}.$$

The expected number of vertices of $\Pi_{N + 1}$ is equal to the expected
number of edges of $\Pi_{N + 1}$. We thus consider the probability that 2
points $P_1, P_2$ chosen from the $N + 1$ random points and their
antipodes define an edge of $\Pi_{N + 1}$. If $P_1 = -P_2$, then they
define an edge with probability 0. Otherwise, the probability that they
define an edge is the probability that the other random points and their
antipodes all lie on the same side of the line $\overline{P_1 P_2}$, which
is the case if the $N - 1$ other random points all lie in the strip
between this line and its reflection in the origin. There are $\binom{2 (N
+ 1)}{2}-(N + 1) = 2 N (N + 1)$ pairs of points which are not antipodal.
Therefore we have
$$V_{N + 1} = 2 N (N + 1) \int_K \int_K A(P_1, P_2)^{N-1} dP_1 dP_2,$$
where $A(P_1, P_2)$ is the area of the intersection of $K$ with the strip
described above.

$A(P_1, P_2)$ depends only on the line $\overline{P_1, P_2}$. If this is
the line $\{(x, y) : (x, y) \cdot (\cos \theta, \sin \theta) = r \}$ for
$r \ge 0$, $0 \le \theta < 2 \pi$, then $A(P_1, P_2) = A(r, \theta)$. Now
$$\bp \cos \theta & \sin \theta \\ -\sin \theta & \cos \theta
    \ep \bp \cos \theta \\ \sin \theta \ep
    = \bp 1 \\ 0 \ep,
$$
so the rotation above takes the line $\overline{P_1 P_2}$ to the vertical
line through $(r,0)$. Now if $P_i = (x_i, y_i)$ for $i = 1, 2$, we denote
$$\bp r \\ s_i \ep = \bp \cos \theta & \sin \theta \\ -\sin \theta &
    \cos \theta \ep \bp x_i \\ y_i \ep,$$
so that
$$\bp x_i \\ y_i \ep = \bp \cos \theta & - \sin \theta \\ \sin \theta &
    \cos \theta \ep \bp r \\ s_i \ep
    = \bp r\cos \theta - s_i \sin \theta \\ r \sin \theta
    + s_i \cos \theta \ep.$$
From this follows
$$dx_1 dy_1 dx_2 dy_2 =
    = |s_1 - s_2| dr d\theta ds_1 ds_2.$$
Since
$$\int_a^b \int_a^b |s_1 - s_2| ds_1 ds_2 = \frac{1}{3}(b - a)^3,$$
we have
$$\int_K \int_K A(P_1, P_2)^{N-1} dP_1 dP_2
    = \frac{1}{3} \int_0^{2 \pi}\int_0^\infty A(r, \theta)^{N-1}
        \ell(r, \theta)^3 dr d\theta,$$
since $\ell(r, \theta)$ is the length of the intersection of the line
$\{(x, y) : (x, y) \cdot (\cos \theta, \sin \theta) = r \}$ with $K$. Note
that there is no need to restrict the domain of the integrals on the right
hand side above, since the integrand is automatically 0 outside the domain
of integration.
\end{proof}

\begin{cor} \label{parallelogram}
$$\E V_{\Par, N} = 1 - \frac{4}{3 (N+1)} \sum_{k=1}^{N+1} \frac{1}{k}$$
for each $N \ge 2$.
\end{cor}
\begin{proof}
By symmetry, the integral over $0 \le \theta \le 2 \pi$ in Proposition
\ref{n=2} is 8 times the integral over $0 \le \theta \le \pi / 4$.
For $0 < \theta < \pi / 4$, we have
\begin{align*}
\ell(r, \theta) & = \sec \theta, \\
A(r, \theta) & = 2r \sec \theta
\end{align*}
for $0 < r < \frac{1}{2} (\cos \theta - \sin \theta)$;
\begin{align*}
\ell (r, \theta) & = \left( frac{1}{2} - \frac{r - \
	frac{1}{2} \cos \theta}{\sin \theta} \right) \sec \theta, \\
A(r, \theta) & = 1 - \left( frac{1}{2} - \frac{r - \
        frac{1}{2} \cos \theta}{\sin \theta} \right)^2 \tan \theta
\end{align*}
for $\frac{1}{2}(\cos \theta - \sin \theta) < r < \frac{1}{2}
(\cos \theta + \sin \theta)$; and $\ell (r, \theta) = 0$ for
$r > \frac{1}{2}(\cos \theta + \sin \theta)$.
Using these, the remainder of the proof is elementary integration.
\end{proof}
Similar expressions for $\E U_{\Par, N}$ and $\E U_{\Delta, N}$ for $N \ge 3$
were derived by Buchta \cite{Buchta}.

\begin{cor} \label{n=2-ball}
$$\E V_{\Ell^2, N} = 1 - \frac{2N}{3\pi^N} \int_0^\pi
    (t + \sin t)^{N - 1} (1 + \cos t)^2 dt$$
for each $N \ge 2$.
\end{cor}
\begin{proof}
We may assume that $\Ell^2$ is the disc of radius $R = \pi^{-1/2}$. Then
$\ell(r, \theta)$ and $A(r, \theta)$ are independent of $\theta$. To apply
Proposition \ref{n=2}, we need to compute
$$\int_0^R A(r)^{N - 1} \ell(r)^3 dr
    =R \int_0^{\pi/2} A(R \sin t)^{N - 1} \ell(R \sin t)^3 \cos t \ dt.$$
Now $\ell(R \sin t) = 2 R \cos t$, and we have
\begin{equation*}
A(R \sin t) = 2\int_0^{R \sin t} \ell(s) ds
    = R^2 (2 t + \sin 2t).
\end{equation*}
The claim now follows from Proposition \ref{n=2}.
\end{proof}

From this we calculate the first few values of $\E V_{\Ell^2, N}$:
\begin{align*}
\E V_{\Ell^2, 2} & = \frac{16}{9 \pi^2} \approx 0.1801, &
    \E V_{\Ell^2, 3} & = \frac{35}{12 \pi^2} \approx 0.2955,\\
\E V_{\Ell^2, 4} & = \frac{-5632 + 1575 \pi^2}{270 \pi^4} \approx 0.3769,
& \E V_{\Ell^2, 6} & = \frac{7(-3289 + 600 \pi^2)}{432 \pi^4} \approx
    0.4380.
\end{align*}
A similar expression for $\E U_{\Ell^2, N}$ for $N \ge 3$ was derived by
Efron \cite{Efron}.

\subsection{Expected volume in an ellipsoid} \label{S:Ellipsoid}

Now let $K \in \K_s^3$. For $r \ge 0$, $0 \le \theta<2\pi$, $0 \le \phi <
\pi$, let
\begin{eqnarray*}
H(r, \theta, \phi) & = & \{(x, y, z) : (x, y, z) \cdot
    (\sin \phi \cos \theta, \sin \phi \sin \theta, \cos \phi) = r \}, \\
A(r, \theta, \phi) & = & \vol_2 (K \cap H(r, \theta, \phi)), \\
V(r, \theta, \phi) & = & \vol_3 \{(x, y, z) : |(x, y, z) \cdot
    (\sin \phi \cos \theta, \sin \phi \sin \theta, \cos \phi)| \le r \} \\
    & = & 2 \int_0^r A(s, \theta, \phi) ds,
\end{eqnarray*}
and let $a(r, \theta, \phi) = \E U_{K \cap H(r, \theta, \phi), 3}$.

\begin{prop}\label{n=3}
Let $K \in \K_s^3$ and $N \ge 3$. If $|K| = 1$ then
\begin{multline*}
\E V_{K, N} = 1 - \frac{1}{N + 1} \\
    - \frac{2 N (N - 1)}{3} \int_0^{2 \pi} \int_0^\pi \int_0^\infty
    V(r, \theta, \phi)^{N - 2} A(r, \theta, \phi)^3 a(r, \theta, \phi)
    \sin \phi \ dr d\phi d\theta.
\end{multline*}
\end{prop}
\begin{proof}
The basic approach is the same as in the two-dimensional case. We consider
a random polyhedron $\Pi_{N + 1}$ in $K$ which is the symmetric convex
hull of $N + 1$ independent random points uniformly distributed in $K$.
Let $v_{N + 1}$, $e_{N + 1}$, $f_{N + 1}$ denote the expected number of
vertices, edges, and faces, respectively, of $\Pi_{N + 1}$. Each of the $N
+ 1$ random points is a vertex of $\Pi_{N + 1}$ iff it is not contained in
the symmetric convex hull of the other $N$ random points, therefore it is
a vertex with probability $1 - \E V_{K, N}$. Therefore
$$V_{N + 1} = 2 (N + 1) (1 - \E V_{K, N}).$$
$\Pi_{N + 1}$ is simplicial with probability 1, which implies $e_{N + 1} =
\frac{3}{2} f_{N + 1}$. Together with Euler's formula $v_{N + 1} - e_{N +
1} + f_{N+1} = 2$, these facts imply
$$\E V_{K, N} = 1 - \frac{1}{N + 1} - \frac{1}{4 (N + 1)} f_{N+1}.$$

Now choose three points $P_1, P_2, P_3$ from the $N + 1$ random points and
their antipodes, such that no two of the chosen points are antipodes.
There are $2^3 \binom{N + 1}{3}$ such possible choices. The points $P_1,
P_2, P_3$ span a face of $\Pi_{N + 1}$ iff all of the other random points
and their antipodes lie in the slab between the plane $H(P_1, P_2, P_3)$
containing $P_1, P_2, P_3$ and its opposite. Therefore
$$f_{N + 1} = 8 \binom{N + 1}{3} \int_K \int_K \int_K
    V(P_1, P_2, P_3)^{N-2} dP_1 dP_2 dP_3,$$
where $V(P_1, P_2, P_3)$ is the volume of the intersection of $K$ with the
slab described above.

$V(P_1, P_2, P_3)$ depends only on the plane $H(P_1, P_2, P_3)$. If
$H(P_1, P_2, P_3) = H(r, \theta, \phi)$, then we change variables by first
rotating by (need geometric description here). This will take $H(P_1, P_2,
P_3)$ to the plane parallel to the $x y$ plane through the point $(r, 0,
0)$, that is, to the plane $H(r, 0, \pi/2)$. If $P_i = (x_i, y_i, z_i)$ is
taken to $(r, s_i, t_i)$ by these rotations for $i = 1, 2, 3$, then we
have
$$\bp x_i \\ y_i \\ z_i \ep =
    \bp \cos \theta & - \sin \theta & 0 \\
        \sin \theta & \cos \theta   & 0 \\
        0           & 0             & 1 \ep
    \bp \sin \phi & 0 & - \cos \phi \\
        0         & 1 & 0           \\
        \cos \phi & 0 & \sin \phi   \ep
    \bp r \\ s_i \\ t_i \ep.$$
This change of variables has the Jacobian
$$\begin{Vmatrix} 1 & 1 & 1 \\ s_1 & s_2 & s_3 \\
    t_1 & t_2 & t_3 \end{Vmatrix} \sin \phi.$$
The claim now follows since
$$\begin{Vmatrix} 1 & 1 & 1 \\ s_1 & s_2 & s_3 \\
    t_1 & t_2 & t_3 \end{Vmatrix} $$
is twice the area of the convex hull of $P_1, P_2, P_3$. As in the proof
of Proposition \ref{n=2}, there is no need to restrict the domain of
integration at this point.
\end{proof}

\begin{cor} \label{n=3-ball}
$$\E V_{\Ell^3,N} = 1 - \frac{1}{N + 1} - \frac{105 N (N - 1)}{2^{N + 5}}
    \int_0^1 (1 - t^2)^4 (3 t - t^3)^{N - 2} dt.$$
\end{cor}
\begin{proof}
We may assume that $\Ell^3$ is the ball of radius $R = (\frac{3}{4
\pi})^{1/3}$. Then for any $\theta, \phi$, and $r < R$, $H(r, \theta,
\phi)$ is a disc of radius $\sqrt{R^2 - r^2}$, so
\begin{eqnarray*}
A(r) & = & \pi(R^2 - r^2), \\
a(r) & = & \frac{35}{48 \pi^2} A(r) = \frac{35}{48 \pi}(R^2 - r^2), \\
V(r) & = & 2\pi \int_0^r (R^2 - s^2) ds
    = 2 \pi \left(R^2 r - \frac{1}{3} r^3 \right).
\end{eqnarray*}
The claim then follows from Proposition \ref{n=3}.
\end{proof}

From this we calculate the first few values of $\E V_{\Ell^3, N}$:
\begin{align*}
\E V_{\Ell^3, 3} & = \dfrac{27}{512}, &
    \E V_{\Ell^3, 4} & = \dfrac{72}{715}, \\
\E V_{\Ell^3, 5} & = \dfrac{585}{4096}, &
     \E V_{\Ell^3, 6} & = \dfrac{58104}{323323}.
\end{align*}
A similar expression for $\E U_{\Ell^3, N}$ for $N \ge 4$ was derived by
Efron \cite{Efron}.

Buchta and Reitzner \cite{BR} use a nonsymmetric analogue of Proposition
\ref{n=3} to derive an expression for $\E U_{T, N}$ for $N \ge 4$, where
$T$ is a tetrahedron. It is natural to ask whether Proposition \ref{n=3}
can be used to calculate $\E V_{K, N}$ when $K$ is a cube or octahedron.
The chief difficulty comes from the appearance of the quantity $a(r,
\theta, \phi)$ in the integrand, which depends in general on the shape of
the planar sections of $K$. In the case of the tetrahedron, these sections
are either triangles or quadrilaterals, for which formulas for the
expected area of the convex hull of three random points are known. For
polyhedra with more facets, planar sections can be polygons for which the
necessary values of $a(r, \theta, \phi)$ are not known.

Unfortunately, it does not seem feasible to extend directly the approach
in this and the previous section to $n \ge 4$. The reason is that the
proofs of Propositions \ref{n=2} and \ref{n=3} actually calculate the
expected number of facets of $\Pi_{N + 1}$, whereas $\E V_{K, N}$ is
directly related to the expected number of vertices of $\Pi_{N + 1}$. In
the plane, these are equal, and in $\R^3$ they are related via Euler's
formula with the fact that $\Pi_{N + 1}$ is almost surely simplicial. If
$n \ge 4$ however, the number of facets of a simplicial polytope does not
uniquely determine the number of vertices.

\section*{Acknowledgements}
This paper is part of the author's Ph.D. thesis, written under the
supervision of Profs. S. Szarek and E. Werner. The author would like to
thank M. Reitzner for suggesting that the techniques of \cite{CCG} should
be useful in the present context.

\end{document}